\documentclass[a4paper,11pt,dvipdfmx]{amsart}
\usepackage[svgnames]{xcolor}% tikzより前に読み込む必要あり
\usepackage{tikz}

\usepackage{comment}

\usepackage[hang,small,bf]{caption}
\usepackage[subrefformat=parens]{subcaption}
\usepackage{etoolbox}
\ifdef{\crop}{%
\usepackage[includeheadfoot,twoside=False,paperwidth=448pt,paperheight=587pt,rmargin=15pt,lmargin=15pt,tmargin=15pt,bmargin=15pt]{geometry}%
}{%
\setlength{\topmargin}{22mm}
\addtolength{\topmargin}{-1in}
\setlength{\oddsidemargin}{27mm}
\addtolength{\oddsidemargin}{-1in}
\setlength{\evensidemargin}{27mm}
\addtolength{\evensidemargin}{-1in}
\setlength{\textwidth}{156mm}
\setlength{\textheight}{240mm}
}%

\usepackage{mathtools}
\usepackage{color}
\usepackage[active]{srcltx}
\usepackage{amsmath,amsthm,amsxtra}
\usepackage{amssymb}

\usepackage[mathscr]{eucal}

\usepackage{bm}
\usepackage[all]{xy}

\usepackage{aliascnt}

\usepackage{tabularx}
\newcolumntype{C}{>{\centering\arraybackslash}X} %中央揃え

\theoremstyle{plain}
\newtheorem{thm}{Theorem}[section]
\newtheorem*{thm*}{Theorem}

\newaliascnt{prop}{thm}
\newaliascnt{cor}{thm}
\newaliascnt{lem}{thm}
\newaliascnt{claim}{thm}
\newaliascnt{defn}{thm}
\newaliascnt{ques}{thm}
\newaliascnt{conj}{thm}
\newaliascnt{fact}{thm}
\newaliascnt{rem}{thm}
\newaliascnt{ex}{thm}
\newaliascnt{sett}{thm}

\newtheorem{cor}[cor]{Corollary}

\newtheorem*{prop*}{Proposition}
\newtheorem*{cor*}{Corollary}
\newtheorem*{lem*}{Lemma}
\newtheorem*{claim*}{Claim}
\theoremstyle{definition}

\newtheorem*{defn*}{Definition}
\newtheorem*{ques*}{Question}
\newtheorem*{conj*}{Conjecture}

\newtheorem*{prob*}{Problem}
\newtheorem{rem}[rem]{Remark}
\newtheorem{ex}[ex]{Example}
\newtheorem*{fact*}{Fact}
\newtheorem*{rem*}{Remark}
\newtheorem*{ex*}{Example}
\aliascntresetthe{prop}
\aliascntresetthe{cor}
\aliascntresetthe{lem}
\aliascntresetthe{claim}
\aliascntresetthe{defn}
\aliascntresetthe{ques}
\aliascntresetthe{conj}
\aliascntresetthe{fact}
\aliascntresetthe{rem}
\aliascntresetthe{ex}
\aliascntresetthe{sett}

\usepackage[pointedenum]{paralist}
\usepackage{varioref}
\labelformat{equation}{\textnormal{(#1)}}
\labelformat{enumi}{\textnormal{(#1)}}

\usepackage[a4paper]{hyperref}

\def\textsectionN~{\textsection{}}

\renewcommand\phi{\varphi}
\renewcommand\epsilon{\varepsilon}
\renewcommand\leq{\leqslant}
\renewcommand\geq{\geqslant}

\makeatletter
\newcommand{\set}{  \@ifstar{\@setstar}{\@set}}\newcommand{\@setstar}[2]{\{\, #1 \, ; \,  #2 \,\}}
\newcommand{\@set}[1]{\{ #1 \}}
\makeatother

\newcommand{\trans}[1][1]{\raisebox{#1ex}{\scriptsize\kern0.1em$t$\kern-0.1em}}

\DeclareMathOperator{\Hom}{Hom}

\newcommand{\wid}{\mathrm{width}}

\DeclareMathOperator{\mult}{mult}

\DeclareMathOperator{\vol}{vol}
\DeclareMathOperator{\Conv}{Conv}

\DeclareMathOperator{\GL}{GL}

\def\Z{\mathbb{Z}}
\def\Q{\mathbb{Q}}
\def\R{\mathbb{R}}
\def\C{\mathbb{C}}

\def\r+{\mathbb{R}_{\geq 0}}

\def\ep{\varepsilon}

\def\r+{{\R}_{\geq 0}}
\def\q+{{\Q}_{\geq 0}}
\def\P{\mathbb{P}}

\def\*c{\C^{\times}}

\def\<{\langle}
\def\>{\rangle}

\def\C{\mathbb {C}}

\def\Q{\mathbb {Q}}
\def\R{\mathbb {R}}

\def\Z{\mathbb {Z}}

\newcommand{\calo}{\mathcal {O}}

\makeatletter
  
  \@addtoreset{equation}{section}
\makeatother

\title
[A note on Seshadri constants, Gromov widths and lattice widths]
{A note on Seshadri constants, Gromov widths of toric surfaces and lattice widths of polygons}

\author[A.~Ito]{Atsushi~Ito}
\address{Department of Mathematics, 
Faculty of %Natural Science and Technology,
Environmental, Life, Natural Science and Technology,
Okayama University,
Okayama, Japan}
\email{ito-atsushi@okayama-u.ac.jp}

\begin{document}

\maketitle

\begin{abstract}
In this note, 
we study Seshadri constants and Gromov widths of toric surfaces via lattice widths of their moment polygons.
We give the sharp lower bound of the ratio between the Gromov width of a symplectic toric $4$-fold and the lattice width of the moment polygon,
which answers to a question raised by Codenotti, Hall and Hofscheier.
%in  \cite{Codenotti:2021ub}.
\end{abstract}

\section{Introduction}

\vspace{2mm}

For a polarized toric surface $(X,L)$ or a $4$-dimensional symplectic toric manifold $(X,\omega)$, 
we can define the moment polygon $P \subset \R^2$.
In this note, we consider the following three invariants 
for toric surfaces: 
\begin{itemize}
\item the Seshadri constant $\ep(X,L;p)$ for $p \in X$ in algebraic geometry, 
\item the Gromov width $ w_G(X,\omega) $ in symplectic geometry, 
\item the lattice width $ \wid(P)$ in convex geometry.
\end{itemize}
See \S\ref{sec_pre} for the definitions of these invariants.
We just note here that it is difficult to compute 
Seshadri constants and Gromov widths in general, even for toric surfaces.
On the other hand,  it is easier to compute lattice widths.

For a $4$-dimensional symplectic manifold $(X,\omega) $ with moment polygon $P \subset \R^2$,
we have 
\begin{align*}
\frac12 \wid(P) \leq w_G(X,\omega) \leq \wid(P)
\end{align*}
by \cite{Codenotti:2021ub} and  \cite{Chaidez:2020aa}.
In \cite[page 3]{Codenotti:2021ub}, it is asked what the sharp lower bound of the ratio $w_G(X_P, \omega_P)/ \mathrm{width}(P)$.
\footnote{In fact, we have $\frac23 \wid(P) \leq w_G(X,\omega) $ by \cite[Sections 5, 6]{MR4065714}. }

To answer to this question,
we consider the corresponding question for Seshadri constants instead of Gromov widths
since $\ep(X_P,L_P;p) \leq w_G(X,\omega) $ holds for the toric surface $X_P$, the ample $\R$-divisor $L_P$ on $X_P$
corresponding to $P$ and any $p \in X_P$.
For a convex body $\Delta \subset \R^2$, that is, a compact convex set with non-empty interior,
an invariant $\ep(\Delta;1) >0$ is defined in \cite[Section 4]{MR3053712} so that the following are satisfied:
\begin{enumerate}[(i)]
\setlength{\itemsep}{0mm}
\item  If $\Delta$ is a rational polygon, $\ep(\Delta;1) =\ep(X_{\Delta}, L_{\Delta} ;1_{\Delta}) $ holds 
for the $\Q$-polarized toric surface $ (X_{\Delta}, L_{\Delta} )$ corresponding to $\Delta$
and the identity $1_{\Delta} \in (\C^*)^2 \subset X_{\Delta}$ of the maximal torus.
\item For a convex body $\Delta \subset \R^2$ and an increasing sequence of convex bodies 
$ \Delta_1 \subset \Delta_2 \subset \cdots \subset \Delta_i \subset \cdots \subset \R^2$ with $\Delta = \overline{\bigcup_{i=1}^{\infty} \Delta_i} $,
it holds that $\ep(\Delta;1) =\sup_{i } \ep(\Delta_i;1) = \lim_i \ep(\Delta_i;1) $.
\end{enumerate}

The purpose of this note is to show the following theorem and corollary.

\begin{thm}\label{main_thm}
Let $\Delta \subset  \R^2$ be a convex body.
Then
\begin{enumerate}
\item $\frac34 \, \mathrm{width}(\Delta) \leq \ep(\Delta;1) \leq  \mathrm{width}(\Delta)$ holds.
\item 
The equality $ \ep(\Delta;1)= \frac34 \, \mathrm{width}(\Delta)$ holds if and only if, up to unimodular transformations, 
$\Delta $ is a parallel translation of $ t P_0$ for some $ t >0$, 
where $P_0 \subset \R^2$ is the convex hull of $\{(1,0), (0,1), (-1,-1)\}$.
\end{enumerate}
\end{thm}

We note that $(X_{P_0}, L_{P_0}) $ coincides with $(S,\calo_S(1))$, where $S$ is the cubic surface $S \subset \P^3$ defined by $xyz-w^3=0$.

\begin{cor}\label{cor_GW}
Let $(X,\omega)$ be a $4$-dimensional  symplectic toric manifold with moment polygon 
$P \subset \R^2$.
Then 
\begin{enumerate}
\item $\frac34 \,  \mathrm{width}(P) < w_G(X, \omega) \leq   \mathrm{width}(P)$ holds.
\item The lower bound in (1)  is sharp in the following sense;
for any $\ep >0$, there exists $(X,\omega)$ such that $w_G(X, \omega)/ \mathrm{width}(P) < \frac34 + \ep$.
\end{enumerate}
\end{cor}

The reason why the bound $ \frac34 \mathrm{width}(P) $ is not attained in \autoref{cor_GW} (1)  is that 
$P_0$ in \autoref{main_thm} (2) is not Delzant.

In this note, varieties are defined over the field of complex numbers $\C$.

\subsection*{Acknowledgments}
The author would like to thank Professor Florin Ambro for useful comments.
The author was supported by JSPS KAKENHI Grant Numbers 17K14162, 21K03201.

\section{Preliminaries}\label{sec_pre}

First, we recall some notation of toric geometry. We refer the readers to \cite{MR2810322}, \cite{MR1234037} for the detail.

Let $M \simeq \Z^n$ be a free abelian group and  $P \subset M_{\R}$ be a lattice polytope of dimension $n$.
Let $X_P$ be the toric variety defined by the normal fan $\Sigma_P$ of $P$.
Let $v_1,\dots, v_r$ be all the primitive vectors in the dual lattice $N=\Hom(M,\Z)$ which span rays in $\Sigma_P$.
Then we can write
\[
P=\{ u \in M_{\R} \, ; \,   \langle v_i, u \rangle \geq -a_i \text{ for } 1 \leq i \leq r\}
\]
for some $a_i \in \Z$ and we define a divisor $L_P$ on $X_P$ by
\[
L_P= \sum_{i=1}^r a_i D_i,
\]
where $D_i \subset X_P$ is the prime divisor corresponding to the ray spanned by $v_i$.
It is known that $L_P$ is an ample Cartier divisor on $X_P$.

This construction can be generalized to the case when $P $ is a \emph{quasi-rational polytope},
that is,  a polytope whose normal fan is rational.
In fact, the normal fan $\Sigma_P$ is rational as well and hence we can define a toric variety $X_P$ and an $\R$-divisor $L_P$ in a similar manner.
In \cite[Section 4]{MR3466351}, it is shown that $L_P$ is an ample $\R$-Cartier $\R$-divisor.

\subsection{Seshadri constants}

The \emph{Seshadri constant} of a polarized variety $(X, L)$ at a point $p \in X $ is defined as
\begin{align}\label{eq_def_SC}
 \ep(X,L;p)=\inf_C \left\{\dfrac{(C.L)}{\mult_p(C)}\right \},
\end{align}
where we take the infimum for all reduced and irreducible curves $C$ on $X$ passing through $p$, $(C.L) $ is the intersection number, and $\mult_p(C)$ is the multiplicity of $C$ at $p$.
This invariant measures the local positivity of $L$ at $p$.
We note that we can define $ \ep(X,L;p)$ in the same way
when $(X,L)$ is an \emph{$\R$-polarized variety},
that is, when $L$ is an ample $\R$-Cartier $\R$-divisor on $X$.

We can define $\ep(\Delta; 1)  > 0$ for a convex body $\Delta \subset M_{\R}$ 
so that (i), (ii) in the introduction are satisfied by \cite[Section 4]{MR3053712}.
Approximating $\Delta$ by rational polytopes
and applying (ii),
we see that $\ep(\Delta;1) =\ep(X_{\Delta}, L_{\Delta} ;1_{\Delta}) $ in (i) still holds when  $\Delta \subset M_{\R}$ is  a quasi-rational polytope.

\subsection{Gromov widths}

The \emph{Gromov width} $w_G(X, \omega)$ of a symplectic manifold $(X,\omega)$ of dimension $2n$ is the supremum of $\pi r^2$ 
such that the open ball $B(r) =\{ x \in \R^{2n} \, ; \,  ||x|| < r\} $ with the standard symplectic form $\omega_{\mathrm{std}} = \sum_{i=1}^n dx_{2i-1} \wedge dx_{2i}$ can be  symplectically embedded in $X$.
For a smooth $\R$-polarized variety $(X,L)$, let $\omega_L$ be a K\"{a}hler  form on $X$ representing the chern class $c_1(L) \in H^2(X,\R)$.
Then 
\begin{align}\label{eq_SC_leq_GW}
\ep(X,L;p) \leq w_G(X,\omega_L)
\end{align}
holds for any $p \in X$ by \cite[Corollary 2.1.D]{MR1262938}  (see also \cite[Theorem 5.1.22, Proposition 5.3.17]{MR2095471}).

\subsection{Lattice widths}\label{subsec_lattice_width}

For a convex body $\Delta \subset M_{\R}$, 
the \emph{lattice width} $\wid(\Delta) $  is defined as
\[
\wid(\Delta) =  \min_{v \in N \setminus \{0\} }\max_{u,u' \in \Delta} | v(u) -v(u') | ,
\]
where we regard $v \in N \subset N_{\R} $ as a linear map $v : M_{\R} \to \R$. 
Equivalently, $\wid(\Delta)$ is the minimum of the length of the image  $v(\Delta) \subset \R$ for $v \in N \setminus \{0\}$.

From the definitions of Seshadri constants and lattice widths,
we see that $\ep(\Delta;1) \leq \wid(\Delta)$ holds.\footnote{For instance, see \cite[Theorem 3.6]{MR3276156} when $\Delta$ is a rational polytope.
General case follows by approximating the convex body $\Delta$ by rational polytopes.}

For a $2n$-dimensional symplectic manifold $(X,\omega) $ with moment polytope $\Delta \subset M_{\R}$, 
the corresponding inequality
$w_G(X,\omega) \leq \wid(\Delta)$ is proved by \cite{Chaidez:2020aa} for $n =2$ and by  \cite{CHARY2024105149} for general $n \geq 1$.

\section{Proof of \autoref{main_thm}}\label{sec_proof_of_Main_thm}

\begin{proof}[Proof of \autoref{main_thm}]

Let $M=\Z^2$.
For convex bodies $\Delta, \Delta' \subset  M_{\R}=\R^2$,
we say that $\Delta$ is \emph{equivalent} to $ \Delta'$, and write $\Delta \sim  \Delta'$, 
if there exist $g \in \GL_2(\Z)$ and $ u\in \R^2$ such that $\Delta=g( \Delta') +u $.
By  \cite[Section 4]{MR3053712},  $\ep(\Delta;1) =\ep(\Delta';1)$ holds if $\Delta \sim  \Delta'$.

Let $\Delta \subset  \R^2$ be a convex body and set $w =\wid(\Delta)$.
As stated in \autoref{subsec_lattice_width},
$\ep(\Delta;1) \leq w$ holds.
Hence the rest is to show $\ep(\Delta;1) \geq \frac34w$ with equality if and only if $\Delta \sim t P_0$ for some $t>0$,
where $P_0=\Conv ((1,0), (0,1), (-1,-1)) \subset \R^2$.
It is easy to check $\wid(P_0)=2$ and hence
$\Delta \sim t P_0$ implies $ t= \frac{w}2$.

It is known that $ \ep(P_0;1) = \frac32$ by \cite[p.153]{MR3276156}. 
Hence, if $\Delta \sim \frac{w}{2} P_0$,  we have $\ep(\Delta;1)=\ep(\frac{w}{2}  P_0 ; 1) =\frac{w}{2}  \ep(P_0 ; 1)= \frac34 w$. 
Thus it suffices to show $\ep(\Delta;1) > \frac34 w $ if $\Delta \not \sim \frac{w}{2} P_0$.

Assume $\Delta \not \sim \frac{w}{2} P_0$.
Then $ w < \sqrt{\frac{8}{3} \vol(\Delta)} $ 
by \cite[Theorem 2]{MR370369} (see also \cite[Theorem 2.3]{MR2890359}),
where $ \vol (\Delta)$ is the Euclidean volume of $\Delta \subset \R^2$.
Hence if $\ep(\Delta;1) \leq  \frac34 w $,
we have $\ep(\Delta;1) < \sqrt{\frac{3}{2} \vol(\Delta)} $.
Let $P \subset \Delta$ be a rational polytope such that $P$ sufficiently approximates $\Delta$.
Then $ \ep(P;1),  \wid(P)$, and $\vol(P)$ are sufficiently close to $ \ep(\Delta;1),  w=\wid(\Delta)$, and $\vol(\Delta)$, respectively.
Since $ \ep(\Delta;1) \leq \frac34 w <w$, 
it holds that  $\ep(P;1) <  \wid(P)$ and $\ep(P;1) <\sqrt{\frac{3}{2} \vol(P)} =\sqrt{\frac{3}{4} (L_{P}^2)}$,
where $(L_P^2)$ is the self intersection number of $L_P$.
We can take such $P$ so that $X_{P}$ is smooth.
Since $ \ep(P;1) =\ep(X_{P}, L_{P}; 1_{P}) = \ep(X_{P}, L_{P}; p)  $ holds for any point $p $ in the maximal torus of $X_{P}$, 
we can apply \cite[Theorem 2]{MR1962050} to $(X_{P}, L_{P})$ 
and obtain a fibration $f : X_{P} \to C$ to a smooth curve such that a general fiber $F$ of $f$ satisfies $ (F.L_{P}) < \wid(P)$.
Since $X_{P}$ is toric, the fibration $f$ is a toric  morphism.
Let $\varpi : N \to \Z$ be the corresponding group homomorphism.
Then we have  the natural surjection $\pi : M \to (\ker \varpi)^{\vee} \simeq \Z$.
A general fiber $F$ of $ f$ is isomorphic to $X_{\pi(P)}$ with $L_{P} |_{F} \equiv L_{\pi(P)}$ and hence $ (F.L_{P}) = |\pi(P)|$ holds.
Thus we have
\begin{align*}
\wid(P) >  (F.L_{P}) = |\pi(P)| \geq \wid(P),
\end{align*}
which is a contradiction.
Hence $\ep(\Delta;1) > \frac34 w$ holds.
\end{proof}

\begin{proof}[Proof of \autoref{cor_GW}]
(1)
The lower bound $\frac34 \wid(P) < w_G(X_P,\omega_P)$
 follows from \autoref{main_thm}  and 
  \ref{eq_SC_leq_GW}. 
We note that $P_0$ in \autoref{main_thm} (2)  is not Delzant and hence
the moment polygon $P$ cannot be equivalent to  a multiple of $P_0$.
The upper bound $w_G(X_P,\omega_P) \leq \wid(P) $ is proved in \cite[Theorem 1.11]{Chaidez:2020aa}. 
(2) follows from \autoref{eg_GW} below.
\end{proof}

\begin{ex}\label{eg_GW}
Let $ P_0=\Conv \{(1,0), (0,1), (-1,-1)\}$ as above.
For an integer $k \geq 1$, let
\begin{align*}
Q_0&=\Conv ( (2,0), (2,1), (1,2), (0,2), (-1,1), (-2,-1), (-2,-2), (-1,-2), (1,-1) ) \\[1mm]
Q_k&=kP_0 + Q_0 \\
&=\Conv ( (k+2,0), (k+2,1), (1,k+2), (0,k+2), (-1,k+1) \\
& \hspace{20mm}(-k-2,-k-1), (-k-2,-k-2), (-k-1,-k-2), (k+1,-1) ).
\end{align*}
\begin{figure}[htbp]
%\vspace{-20mm}
 \begin{minipage}[b]{0.2\linewidth}
  \centering
  \includegraphics[keepaspectratio, scale=0.3]
  {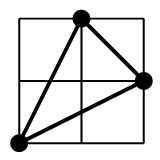}
  \subcaption*{$P_0$}\label{}
 \end{minipage}
 \begin{minipage}[b]{0.32\linewidth}
  \centering
  \includegraphics[keepaspectratio, scale=0.3]
  {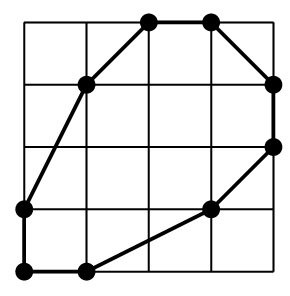}
  \subcaption*{$Q_0$}\label{}
 \end{minipage}
 \begin{minipage}[b]{0.32\linewidth}
  \centering
  \includegraphics[keepaspectratio, scale=0.3]
  {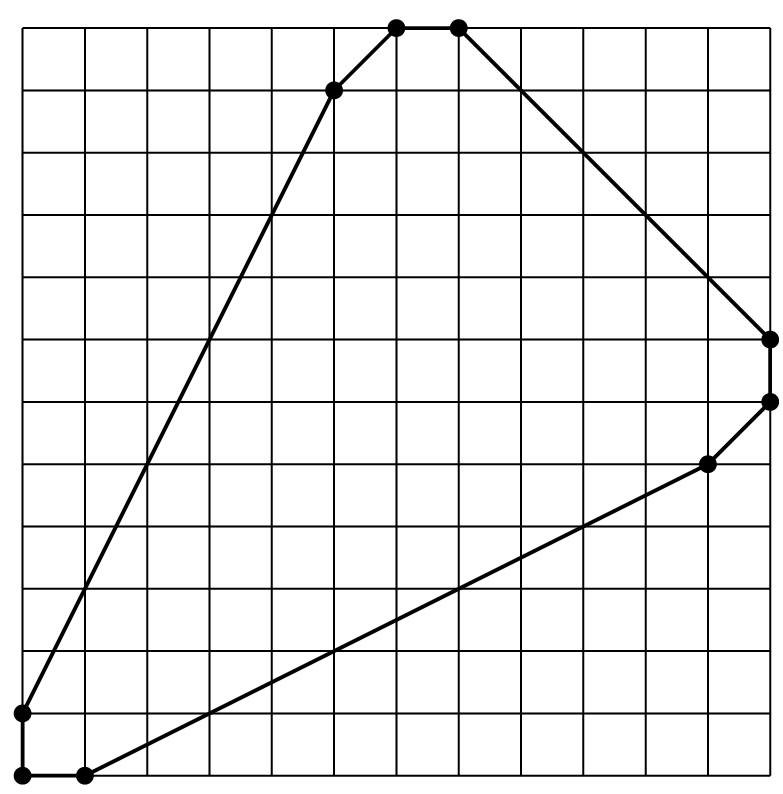}
  \subcaption*{$Q_4$}\label{}
 \end{minipage}
 \caption{}\label{figure_Q_i}
\end{figure}

In \cite[Example 5.6]{MR4251298}, the authors consider $Q_4$ and show $w_G(X_{Q_4},\omega_{Q_4}) =\frac{21}2$ 
using the algorithm from \cite[Section 6]{MR3734610}.
We can compute $w_G(X_{Q_k},\omega_{Q_k}) =\frac{3k+9 }2$ for any $k \geq 1 $ similarly.
On the other hand, we can check $\wid(Q_k)= 2k+4$.
Hence we have
\[
\frac{w_G(X_{Q_k},\omega_{Q_k})}{\wid(Q_k)} =\frac{3k+9}{4k+8} \xrightarrow{k \to \infty } \frac34.
\]
\end{ex}

\begin{rem}
We can also show $\ep(X_{Q_k}, L_{Q_k};1_{Q_k} ) = \frac{3k+9}{2}$ for $k \geq 1$ 
as follows:

The linear system $|L_{P_0}|$ embeds $X_{P_0}$ into $\P^3$ as the cubic surface defined by $xyz-w^3=0$.
Let $H \subset \P^3$ be the hyperplane defined by $x+y+z -3w=0$.
Since $1_{P_0} =[1:1:1:1] \in X_{P_0} \subset \P^3$,
we can check that $C_0\coloneqq  H|_{X_{P_0}} \in |L_{P_0}|$ is an irreducible curve with $\mult_{1_{P_0}} (C_0)=2. $

Since the normal fan of $Q_{0}$ refines that of $P_0$, we have a morphism $\pi : X_{Q_0} \to X_{P_0}$.
Since $ \pi $ is an isomorphism over $C_0$,
the pullback $\tilde{C}_0\coloneqq \pi^* C_0 \subset  X_{Q_0}$ is irreducible.
Since $ (\pi^* L_{P_0}^2 ) =3, (\pi^*L_{P_0}. L_{Q_0})=9$ and $(X_{Q_k}, L_{Q_k}) = (X_{Q_0}, k \pi^* L_{P_0} + L_{Q_0})$,
we have
\begin{align}\label{eq_3k+9}
 \frac{(\tilde{C}_0 . L_{Q_k})}{\mult_{1_{Q_k}} (\tilde{C}_0)} =  \frac{ k ( \pi^* L_{P_0}^2)+ ( \pi^* L_{P_0}  . L_{Q_0})}{2} = \frac{3k+9}2.
\end{align}

Let $C \subset X_{Q_k}$ be an irreducible curve passing through $1_{Q_k}$ with $ C \neq \tilde{C}_0$.
Since $  Q_0$ contains $2 P_0$,
the torus invariant divisor $L_{Q_0} -2 \pi^* L_{P_0} $ is 
effective.
Since $ C$ is not contained in boundary divisors of $X_{Q_0}=X_{Q_k}$, we have $(C.L_{Q_0}) \geq 2 (C.  \pi^* L_{P_0}) $ and hence
\begin{align}\label{eq_2k+4}
\begin{split}
 \frac{(C . L_{Q_k})}{\mult_{1_{Q_k}}( C)} &= \frac{(C . k \pi^* L_{P_0} + L_{Q_0})}{\mult_{1_{Q_k}}( C)} \\
 &\geq  (k+2) \frac{ (C .  \pi^* L_{P_0})}{\mult_{1_{Q_k}}( C)} \\
&= (k+2) \frac{ (C .  \tilde{C}_0)}{\mult_{1_{Q_k}}( C)} \\
&\geq(k+2) \frac{\mult_{1_{Q_k}}( C) \cdot \mult_{1_{Q_k}}( \tilde{C}_0)}{\mult_{1_{Q_k}}( C) }  =2(k+2).
 \end{split}
\end{align}
Since $ 2(k+2) \geq \frac{3k+9}{2}$ for $ k\geq 1$,
we have $\ep(X_{Q_k}, L_{Q_k};1_{Q_k} ) = \frac{3k+9}{2}$ by \ref{eq_3k+9}, \ref{eq_2k+4}.
\end{rem}

\bibliographystyle{amsalpha}
%\bibliography{mainbibs}
\providecommand{\bysame}{\leavevmode\hbox to3em{\hrulefill}\thinspace}
\providecommand{\MR}{\relax\ifhmode\unskip\space\fi MR }
\providecommand{\MRhref}[2]{%
  \href{http://www.ams.org/mathscinet-getitem?mr=#1}{#2}
}
\providecommand{\href}[2]{#2}

\end{document}